\documentstyle[amsfonts,11pt,leqno]{article}
\begin{document}
\title{ Quasi $Q_n$-filiform Lie algebras
\thanks{ Supported by National Natural Science Foundation of China
(grant No. 10571091,10671027) E-mail address: renbin1964@163.com }}
 \author{ Bin Ren $^{1}$  \qquad Lin Sheng Zhu $^2$
 \\ {\small {}$^1$ Department of Mathematics,
 University of Science and Technology of Suzhou,}\\
 {\small Suzhou 215009, China}
  \\ {\small {}$^2$ Department of Mathematics, Changshu college,}\\
 {\small Changshu 215500, China}   }
\date{}

\maketitle
\begin{abstract}  In this paper  we explicitly determine the derivation algebra,
automorphism group of quasi $Q_n$-filiform Lie algebras, and
applying some properties of root vector decomposition we obtain
their isomorphism theorem.


\end{abstract}
\maketitle

Solvable Lie algebras are very important in the study of Lie
algebras by Levi decomposition theorem. During the last decades,
they have also shown their importance in physics. But the
classification of solvable Lie algebras seems to be feasible because
of the absence of their global structural properties. So one focuses
on some specific classes of solvable case. These specific classes
are often related to specific classes of nilpotent Lie algebras,
such as Heisenberg algebras, filiform Lie algebras and so on.

 Filiform Lie algebras is an important class of
nilpotent Lie algebras. In recent  years, this class of  algebras
was studied in many articles. In [11] we study a class of nilpotent
Lie algebras which is related to filiform Lie algebras, we call them
quasi $L_n$-filiform Lie algebras. In [11] we obtain some of their
structural properties. The natural idea is to study the more
nilpotent Lie algebras which are related to filiform Lie algebras.
In this paper we study  quasi $Q_n$-filiform Lie algebras, we
explicitly determine their derivation algebra, automorphism group,
and obtain their isomorphism theorem.

In this paper, all  nilpotent Lie algebras discussed  are finite
dimensional complex nilpotent Lie algebras. We denote the central
descending sequence by  $c^0N=N$, $c^iN=[N,c^{i-1}N]$.

\section {  Preliminaries}

In this section we recall some elementary facts about  nilpotent Lie
algebras and quasi $L$-filiform Lie algebras.

                       \vspace{6pt}

 \noindent {\bf Lemma 1.1} [12]\quad If $N$ is a nilpotent Lie algebra, the two
   following assertions  are equivalent:

   (1) $\{x_1,x_2,\ldots, x_n\}$ is a minimal system of generators;

  (2) $\{x_1+c^1N, x_2+c^1N,\ldots, x_n+c^1N\}$ is a basis for the vector
   space $N/c^1N$.

                       \vspace{6pt}

 If $H$ is a maximal torus on $N$ (a maximal abelian subalgebra of Der$N$
  which consists of semi-simple linear transformations), then $N$ can be
 decomposed into a direct  sum of root spaces with respect  to $H$:
  $ N=\sum_{\alpha \in H^*}N_{\alpha} $.   The scalar
  mult$(\alpha):=$dim $N_{\alpha}$ is called the multiplicity of the root
  $\alpha$. We also  denote $\dim [x]=\dim N_\alpha$ if $x\in N_\alpha$.

                                      \vspace{6pt}

 \noindent {\bf Definition 1.1} [12]\quad   Let $H$ be a maximal
  torus on $N$.  One calls $H$-msg a minimal system of generators
   which consists of root vectors for $H$ .

      \vspace{6pt}

 \noindent  {\bf Definition 1.2} [6]\quad A nilpotent Lie algebra $N$ is called
   quasi-cyclic if $N$ has a subspace $U$ such that $N=U\dot+U^2\dot+
   \cdots\dot+U^k$, where $U^i=[U,U^{i-1}] $.

  \vspace{6pt}

 \noindent  {\bf Lemma 1.2 }[9]\quad
 Let $N$ be a quasi-cyclic nilpotent Lie  algebra,
  $\{ x_1, x_2, \ldots, x_n \}$  be an $H_1$-msg of $N$,
  $\{ y_1, y_2, \ldots, y_n \}$ be an $H_2$-msg of $N$,
  then  $\exists \theta\in{\rm Aut}N$ such that
 $$   ( \ { y_1}  \ \ {y_2} \ \ \cdots  \ \  {y_n} \ )^t=
  A( \ {  \theta(x_1)}  \  \ { \theta(x_2)} \ \ \cdots \ \ { \theta(x_n)} \ )^t, $$
  where $(  y_1 \  y_2 \  \cdots   \ y_n  )^t $
is the transpose of the matrix  $(  y_1 \  y_2 \  \cdots   \ y_n )$,
$A$ is a $n\times n$ invertible matrix.  In particular, if $\dim
[x_i]=1$, $\forall \ i$, then $A$ is  a  monomial matrix (i.e., each
row and each column has  exactly one   nonzero entry).

   \vspace{6pt}

\noindent {\bf Definition 1.3 } [13]\quad  Let $L$ be a
$n$-dimensional  Lie  algebra,  $L$ is called a filiform Lie algebra
if $ \dim c^i L =n-i-1$ for $ 1\le i\le n-1.$

    \vspace{6pt}

 \noindent {\bf Definition 1.4 }[11]\quad Let $L$ be a  $(n+1)$-dimensional
 filiform Lie  algebra, $N$ is called a quasi $L$-filiform
 Lie algebra  if         $N=N_1+N_2+\cdots+N_m $,
  where $N_i\cong L$  $(1\le i \le m)$,  and
  $ [N_{i}, N_{j}]=0$ $(i\not=j)$.
We denote $N$ by $N(L,m)$ or $N(L,m,r),$  $r=\dim c^{n-1}N.$

   \vspace{6pt}

 \noindent {\bf Remark 1.1}\quad The sum in the  decomposition $N=\sum_1^mN_i$
 is not necessarily  direct, so the subalgebras in the  decomposition can have
 nontrivial intersection.

                       \vspace{6pt}

It is easy to know that $N(L,m,r)$ also admits the following
decomposition:
 $$N(L,m,r)=N_1(L,m_1,1)+N_2(L,m_2,1)+\cdots+N_q(L,m_q,1),$$
 where  $N_i\cap N_j=[N_i, N_j]=0$ $( i\ne j )$, $\sum_1^qm_i=m$.
This decomposition is called a $L$-filiform decomposition.

    \vspace{6pt}

 Filiform Lie algebra $Q_n (n=2d+1>3)$ is a $(n+1)$-dimensional  Lie algebra defined in
the basis $\{  x_{0}, x_{1}, \ldots, x_{n} \}$ by
 $$[x_{0}, x_{i}]=x_{i+1}, \quad  [x_{i}, x_{n-i}]=(-1)^{i}x_{n}, \quad i=1,2, \ldots,n-1,$$
 the undefined brackets being zero or obtained by antisymmetry. Obviously $Q_n$ has  another
 basis $\{  e_{0}, e_{1}, \ldots, e_{n} \}$ given by $  e_{0}=x_{0}+x_{1}$, $  e_{i}=x_{i}$, $ i=1,\ldots, n$.
 This basis satisfies
 $$[e_{0}, e_{i}]=e_{i+1}, \quad i=1,2, \ldots,n-2,    $$ $$
 [e_{i}, e_{n-i}]=(-1)^{i}e_{n}, \quad i=1,2, \ldots,n-1.$$
In following what we choose  this basis in discuss.

 Obviously $Q_n$ is a quasi-cyclic nilpotent Lie  algebra. In general, we have  the following lemma.

                               \vspace{6pt}

 \noindent {\bf Lemma 1.3}\quad  $N(Q_n,m)$ is a quasi-cyclic nilpotent Lie algebra.

         \vspace{6pt}

  \noindent {\bf Proof.}\quad Let   $\{ e_{i0}, e_{i1} \}$ be a minimal system of
 generators of $N_i,$ $U$ be the  vector space spanned by
  $\{ e_{10}, e_{11}, \ldots,  e_{m0}, e_{m1} \}$.
  Obviously  $N=\sum_0^{n-1}c^iU.$

  For $x \in c^{k}U \cap  \sum_0^{k-1}c^iU ,$ write
  $x=\sum_1^m( a_{i0}e_{i0}+ a_{i1}e_{i1}+\cdots+a_{i,k-1}e_{i,k-1} ).$
  As $x \in c^{k}U$, we have
  $$ 0=[e_{s1},\underbrace{[ e_{s0}, \cdots, [ e_{s0}}_{n-3},  [ e_{s1}, x ]]
  =a_{s0}e_{sn},  $$
  $$ 0= [e_{s1},\underbrace{ [ e_{s0}, \cdots, [ e_{s0}}_{n-2}, x ]
  =-a_{s1}e_{sn}, $$
  hence $a_{s0}=a_{s1}=0$ for any $ s.$ Then
    $x=\sum_1^m( a_{i2}e_{i2}+\ldots+a_{i,k-1}e_{i,k-1} ).$

  A same argument shows that $a_{sj}=0$ for any $s$, $j$.
  Thus $x=0$. This implies that
   $c^{k}U \cap  \sum_0^{k-1}c^iU=0, \ 1\le k \le  n-1$.  Hence
  $N=U\dot+c^1U\dot+\cdots\dot+c^{n-1}U$.

  \hfill $\Box$

                  \vspace{6pt}

Lemma 1.3 means that  $\{  e_{10}, e_{11}, \ldots, e_{m0},  e_{m1}
\}$ is a minimal system  of generators  of $N(Q_n,m)$. It is easy
to know that $\{ e_{1j},  e_{2j}, \ldots, e_{mj} \}$ is linearly
independent$(1 < j < n)$.  May assume that $\{ e_{10}, e_{11},
\ldots, e_{1,n-1},  \ldots, e_{m0}, e_{m1}, \ldots,  e_{m,n-1},
e_{1n}, \ldots,  e_{rn} \}$
  is a   basis of $N(Q_n,m,r)$, and $e_{in}=\sum_{j=1}^{r} b_{ji}e_{jn}$. Then
  we have
$$\left(  \begin{array}{cccc}
e_{1n}  & e_{2n} & \cdots  & e_{mn}
\end{array}\right)
=\left(  \begin{array}{cccc} e_{1n}  & e_{2n} & \cdots  & e_{rn}
\end{array}\right)
\left(  \begin{array}{cccc} I & B
\end{array}\right). \eqno(1.1)$$

\section  {On  Der$\mathbf{N(Q_n,m, r)}$ }

In this section we explicitly determine the derivation algebra of
$N(Q_n,m, r)$.

            \vspace{6pt}

Let  $\delta$ be a linear transformation of $N(Q_n,m, r)$ such
that
 $$  \delta e_{st}=\sum_{i=1}^{m}\sum_{j=0}^{n-1} c^{st}_{ij}e_{ij}
+\sum_{i=1}^{r}c^{st}_{in}e_{in}, \quad t=0, 1,$$
$$  \delta e_{st}=[\delta e_{s0}, e_{s,t-1}]+[e_{s0}, \delta e_{s,t-1}],
 \quad  2\le t\le n-1,$$
 $$ \delta e_{sn}=-[\delta e_{s1}, e_{s,n-1}]-[e_{s1}, \delta e_{s,n-1}], \quad 1\le s\le r.$$
Then we have
$$\delta e_{st}
=\lambda_{s,t-1}e_{st}+\sum_{j=2}^{n-t}c^{s1}_{sj}e_{s,j+t-1}
+(-1)^tc^{s0}_{s,n-t+1}e_{sn}, \quad  2\le t\le n-1,$$
$$\delta e_{sn}=\lambda_{s}e_{sn}, \quad 1\le s\le r,$$
where  $\lambda_{si}=ic^{s0}_{s0}+c^{s1}_{s1}$,
$\lambda_{s}=(n-2)c^{s0}_{s0}+2c^{s1}_{s1}$.

                        \vspace{6pt}

\noindent  {\bf Theorem 2.1 }\quad  Let $\delta$ be as above,
$e_{in}=\sum_{j=1}^{r} b_{ji}e_{jn}$, then $\delta \in {\rm Der}N$
if and only if   for any  $1\le s, p\le m $,
  $$ \ b_{sj}(\lambda_s-\lambda_j)=0, \quad s > r, \ \ 1\le j\le r,   \eqno (2.1)   $$
 $$\delta e_{s0}=c^{s0}_{s0}e_{s0}+\sum_{i=2}^{n-1}c^{s0}_{si}e_{si}
 +\sum_{i=1}^{r}c^{s0}_{in}e_{in},      \eqno (2.2)     $$
   $$\delta e_{s1}= \sum_{i=1}^{n-2}c^{s1}_{si}e_{si}
 +\sum_{i=1}^{m}c^{s1}_{i,n-1}e_{i,n-1}+\sum_{i=1}^{r}c^{s1}_{in}e_{in},   \eqno (2.3)$$
$$c^{s1}_{si}=0, \quad i=3,5,\ldots,n-2, \eqno (2.4)$$
$$\left(  \begin{array}{cc}
-c^{s1}_{p,n-1}  & c^{p1}_{s,n-1}
\end{array}\right)
\left(  \begin{array}{cccc}
b_{p1}  & b_{p2} & \cdots  & b_{pr}\\
b_{s1}  & b_{s2} & \cdots  & b_{sr}
\end{array}\right)
=0.\eqno (2.5)$$

         \vspace{6pt}

\noindent  {\bf Proof. }\quad  $(\Rightarrow)$: As  $\delta \in
{\rm Der}N$, for any $s>r$, we have
$$\lambda_{s}\sum_{j=1}^{r}b_{js}e_{jn}
 =\lambda_{s}e_{sn}
 =\delta e_{sn}
 =\delta(\sum_{j=1}^{r} b_{js} e_{jn})
 =\sum_{j=1}^{r} b_{js}\lambda_{j} e_{jn}.$$
 Then
$b_{sj}\lambda_s=b_{sj}\lambda_j, \ 1\le j\le r$, i.e., (2.1) holds.

 For any $ \ 1\le s \ne p\le m$, $1<t<n-1$, by
$$0=\delta [e_{s0},e_{s,n-1}]=-c^{s0}_{s1}e_{sn},$$
$$0=\delta [e_{s1}, e_{st}]=c^{s1}_{s0}e_{s,t+1}+((-1)^{n-t}-1)c^{s1}_{s,n-t}e_{sn},$$
$$0=\delta  [e_{s0}, e_{p1}]=c^{s0}_{p0}e_{p2}+c^{s0}_{p,n-1}e_{pn}
+\sum_{j=1}^{n-2}c^{p1}_{sj}e_{s,j+1},$$
$$0=\delta  [e_{s0}, e_{p0}]=-\sum_{j=1}^{n-2}c^{s0}_{pj}e_{p,j+1}
+\sum_{j=1}^{n-2}c^{p0}_{sj}e_{s,j+1},$$
$$0=\delta  [e_{s1}, e_{p1}]=c^{s1}_{p0}e_{p2}-c^{p1}_{s0}e_{s2}
 +\sum_{j=1}^{r}(c^{s1}_{p,n-1} b_{pj}-c^{p1}_{s,n-1}b_{sj})e_{jn},$$
 we have
$c^{s0}_{s1}=c^{s1}_{s0}=c^{s0}_{p0}=c^{s0}_{p,n-1}=c^{s1}_{p0}=c^{p1}_{s0}=0$,
$c^{s1}_{s,n-t}=0$, $t=2,4,\ldots,n-3$,
$c^{p1}_{sj}=c^{s0}_{pj}=c^{p0}_{sj}=0$, $1\le j \le n-2$,
 and
$c^{s1}_{p,n-1}b_{pj}-c^{p1}_{s,n-1}b_{sj}=0$, $1 \le j \le r$.
    Then (2.2), (2.3), (2.4) and (2.5) hold.

          \vspace{6pt}

 $(\Leftarrow)$: We show that $\delta [x, y]=[\delta x, y]+
 [x, \delta y]$ for any $x, y \in N$.  Set
   $$x=\sum_{i=1}^{m} \sum_{j=0}^{n-1} u_{ij}e_{ij}
   +\sum_{i=1}^{r} u_{in}e_{in}, \quad \quad
  y=\sum_{i=1}^{m}\sum_{j=0}^{n-1} v_{ij}e_{ij}
     +\sum_{i=1}^{r} v_{in}e_{in}.  $$

Since $[e_{is}, e_{jt}]=0$ $ (i\ne j)$, we have
  $$\delta [x, y]
=\sum_{i=1}^{m}\sum_{j=1}^{n-2}(u_{i0}v_{ij}-v_{i0}u_{ij})\delta
e_{i,j+1}+\sum_{i=1}^{m}\sum_{j=1}^{n-1}(-1)^ju_{ij}v_{i,n-j}\delta
e_{in}.$$

On the other hand, note that
$$ \sum_{j=2}^{n-1}u_{ij}\delta e_{ij}
=\sum_{j=2}^{n-1}u_{ij}
\left(\lambda_{i,j-1}e_{ij}+\sum_{k=2}^{n-j}c^{i1}_{ik}e_{i,k+j-1}
+(-1)^jc^{i0}_{i,n-j+1}e_{in} \right)$$
$$=\sum_{j=2}^{n-1}u_{ij}\lambda_{i,j-1}e_{ij}+\sum_{s=3}^{n-1}\sum_{t=2}^{s-1}u_{it}c^{i1}_{i,s-t+1}e_{s}
+\sum_{j=2}^{n-1}u_{ij}(-1)^jc^{i0}_{i,n-j+1}e_{in},$$
 we have
 $$ [\delta x, y]=\left[ \sum_{i=1}^{m}u_{i0}\delta e_{i0}, \ y \right]
 + \left[ \sum_{i=1}^{m}u_{i1}\delta e_{i1},\ y \right]
 +\left[ \sum_{i=1}^{m} \sum_{j=2}^{n-1}u_{ij}\delta e_{ij}, \ y \right]$$ $$
=\sum_{i=1}^{m}u_{i0} \left(   \underbrace{
\sum_{j=1}^{n-2}c^{i0}_{i0}v_{ij}e_{i,j+1}}_{A_1}
\underbrace{-\sum_{j=2}^{n-2}c^{i0}_{ij}v_{i0}e_{i,j+1}}_{A_2}
\underbrace{-\sum_{j=1}^{n-2}(-1)^{j}c_{i,n-j}^{i0}v_{ij}e_{in}}_{A_3}
\right)$$
$$+\sum_{i=1}^{m}u_{i1}  \left(
\underbrace{-c^{i1}_{i1}v_{i0}e_{i2}}_{B_1}
\underbrace{-c^{i1}_{i1}v_{i,n-1}e_{in}}_{B_2}
\underbrace{-\sum_{j=2}^{n-2}c^{i1}_{ij}v_{i0}e_{i,j+1}}_{B_3}
\underbrace{-\sum_{j=1}^{n-2}(-1)^{j}c^{i1}_{i,n-j}v_{ij}e_{in}}_{B_4}
\right)$$
$$
+\underbrace{\sum_{i=1}^{m} u_{i1}\sum_{j=1,j\neq i}^{m}
c^{i1}_{j,n-1}v_{j1}e_{jn}}_B$$
$$+\sum_{i=1}^{m} \left(
\underbrace{-\sum_{j=2}^{n-2}u_{ij}\lambda_{i,j-1}v_{i0}e_{i,j+1}}_{E_1}
+\underbrace{\sum_{j=2}^{n-1}(-1)^ju_{ij}\lambda_{i,j-1}v_{i,n-j}e_{in}}_{E_2}
\right)$$
$$+\sum_{i=1}^{m} \left(  \underbrace{\sum_{j=1}^{n-3}(-1)^{j+1}v_{ij}(\sum_{t=2}^{n-j-1}
u_{it}c^{i1}_{i,n-j-t+1})e_{in}}_{E_3}
\underbrace{-\sum_{j=2}^{n-2}u_{ij}\sum_{k=2}^{n-j-1}c^{i1}_{ik}v_{i0}e_{i,k+j}}_{E_4}
\right).$$

  Similarly  we have
$$ [x, \delta y]
=-\sum_{i=1}^{m}v_{i0}  \left( \underbrace{
\sum_{j=1}^{n-2}c^{i0}_{i0}u_{ij}e_{i,j+1}}_{A'_1}
\underbrace{-\sum_{j=2}^{n-2}c^{i0}_{ij}u_{i0}e_{i,j+1}}_{A'_2}
\underbrace{-\sum_{j=1}^{n-2}(-1)^{j}c_{i,n-j}^{i0}u_{ij}e_{in}}_{A'_3}
\right)$$
$$-\sum_{i=1}^{m}v_{i1} \left(
\underbrace{-c^{i1}_{i1}u_{i0}e_{i2}}_{B'_1}
\underbrace{-c^{i1}_{i1}u_{i,n-1}e_{in}}_{B'_2}
\underbrace{-\sum_{j=2}^{n-2}c^{i1}_{ij}u_{i0}e_{i,j+1}}_{B'_3}
\underbrace{-\sum_{j=1}^{n-2}(-1)^{j}c^{i1}_{i,n-j}u_{ij}e_{in}}_{B'_4}
\right)$$
$$
\underbrace{-\sum_{i=1}^{m} v_{i1}\sum_{j=1,j\neq i}^{m}
c^{i1}_{j,n-1}u_{j1}e_{jn}}_{B'}$$
$$-\sum_{i=1}^{m} \left(
\underbrace{-\sum_{j=2}^{n-2}v_{ij}\lambda_{i,j-1}u_{i0}e_{i,j+1}}_{E'_1}
+\underbrace{\sum_{j=2}^{n-1}(-1)^jv_{ij}\lambda_{i,j-1}u_{i,n-j}e_{in}}_{E'_2}
\right)$$
$$-\sum_{i=1}^{m} \left(\underbrace{\sum_{j=1}^{n-3}(-1)^{j+1}u_{ij}(\sum_{t=2}^{n-j-1}
v_{it}c^{i1}_{i,n-j-t+1})e_{in}}_{E'_3}
\underbrace{-\sum_{j=2}^{n-2}v_{ij}\sum_{k=2}^{n-j-1}c^{i1}_{ik}u_{i0}e_{i,k+j}}_{E'_4}
\right).$$

We now  prove that
  $$ [\delta x, y]+[x, \delta y]=B+B'+\delta [x, y]. \eqno (2.6)$$

In order to prove (2.6),  we take
$$\alpha=u_{i1}B_4+E_3-v_{i1}B'_4-E'_3, $$ $$ \beta=u_{i1}B_2-v_{i1}B'_2+E_2-E'_2,  $$
$$\gamma_1=-v_{i1}B'_1+u_{i1}B_1+u_{i0}A_1-v_{i0}A'_1-E'_1+E_1,$$
$$\gamma_2=-v_{i1}B'_3+u_{i1}B_3-E'_4+E_4, $$ $$ \gamma_3=u_{i0}A_3-v_{i0}A'_3.$$

 Since  $\sum_{i=1}^m u_{i0}A_2-\sum_{i=1}^m
v_{i0}A'_2=0$ is obvious, if  we can prove the following
equations: (2.7), (2.8) and (2.9), then  (2.6) holds.
$$\alpha=0, \eqno (2.7)$$
$$ \beta=\sum_{j=1}^{n-1}(-1)^{j}u_{ij}v_{i,n-j}\delta e_{in},\eqno (2.8)
$$ $$
\gamma_1+\gamma_2+\gamma_3=\sum_{j=1}^{n-2}(u_{i0}v_{ij}-v_{i0}u_{ij})\delta
e_{i,j+1}.\eqno (2.9)$$

In fact, note that (2.4), we have
$$\alpha=\sum_{j=1}^{n-2}\sum_{t=1}^{n-j-1}(-1)^{j+1}
u_{it}c^{i1}_{i,n-j-t+1}v_{ij}e_{in}
-\sum_{j=1}^{n-2}\sum_{t=1}^{n-j-1}(-1)^{j+1}
v_{it}c^{i1}_{i,n-j-t+1}u_{ij}e_{in}$$
$$=\sum_{j=1}^{n-2}\sum_{t=1}^{n-j-1}((-1)^{j+1}-(-1)^{t+1})
u_{it}c^{i1}_{i,n-j-t+1}v_{ij}e_{in}=0.$$

Since (2.1) implies that $\delta e_{sn}=\lambda_s e_{sn}$ for any
$s$, and note that
$$E'_2=\sum_{j=2}^{n-1}(-1)^jv_{ij}\lambda_{i,j-1}u_{i,n-j}e_{in}
=-\sum_{j=1}^{n-2}(-1)^{j}v_{i,n-j}\lambda_{i,n-j-1}u_{ij}e_{in},$$
then we have
$$\beta=-c^{i1}_{i1}(u_{i1}v_{i,n-1}-v_{i1}u_{i,n-1} )e_{in}$$ $$
+\sum_{j=2}^{n-1}(-1)^ju_{ij}v_{i,n-j}(\lambda_{i,j-1}+\lambda_{i,n-j-1})e_{in}
$$
$$=\sum_{j=1}^{n-1}(-1)^{j}u_{ij}v_{i,n-j}\lambda_{i}e_{in}.$$

At last it's easy to know that (2.9) holds because of the
following equations:
$$\gamma_1 =(v_{i1}u_{i0}-u_{i1}v_{i0})c^{i1}_{i1}e_{i2}
+\sum_{j=1}^{n-2}(u_{i0}v_{ij}-v_{i0}u_{ij})c^{i0}_{i0}e_{i,j+1}$$
$$+\sum_{j=2}^{n-2}(u_{i0}v_{ij}-v_{i0}u_{ij})\lambda_{i,j-1}e_{i,j+1}
 $$ $$=\sum_{j=1}^{n-2}(u_{i0}v_{ij}-v_{i0}u_{ij})\lambda_{ij}e_{i,j+1},$$
$$\gamma_2
=\sum_{j=2}^{n-2}c^{i1}_{ij}(v_{i1}u_{i0}-u_{i1}v_{i0})e_{i,j+1}
+\sum_{j=2}^{n-2}\sum_{k=2}^{n-j-1}c^{i1}_{ik}(v_{ij}u_{i0}-u_{ij}v_{i0})e_{i,k+j}$$
$$=\sum_{j=1}^{n-2}\sum_{k=2}^{n-j-1}(u_{i0}v_{ij}-v_{i0}u_{ij})c^{i1}_{ik}e_{i,k+j},
$$
$$\gamma_3=\sum_{j=1}^{n-2}u_{i0}(-1)^{j+1}c_{i,n-j}^{i0}v_{ij}e_{in}
-\sum_{j=1}^{n-2}v_{i0}(-1)^{j+1}c_{i,n-j}^{i0}u_{ij}e_{in}$$
$$=\sum_{j=1}^{n-2}(u_{i0}v_{ij}-v_{i0}u_{ij})(-1)^{j+1}c_{i,n-j}^{i0}e_{in}.$$

 Now if we can prove that  $B+B'=0$, then we have
$\delta [x,y]=[\delta x, y]+ [x, \delta y]$.  Next we  prove this
equation. As (2.5) it follows that
  $$2(B+B')=2\sum_{i=1}^{m} \sum_{j=1,j\ne i}^{m}
 ( u_{i1}v_{j1}- v_{i1}u_{j1} )c^{i1}_{j,n-1}e_{jn} $$
$$=\sum_{i=1}^{m} \sum_{j=1,j\ne i}^{m}
 ( u_{i1}v_{j1}- v_{i1}u_{j1} )c^{i1}_{j,n-1}e_{jn}
 +\sum_{i=1}^{m} \sum_{j=1,j\ne i}^{m}
 ( u_{j1}v_{i1}- v_{j1}u_{i1} )c^{j1}_{i,n-1}e_{in}$$
$$=\sum_{i=1}^{m} \sum_{j=1,j\ne i}^{m}
 ( ( u_{i1}v_{j1}- v_{i1}u_{j1} )c^{i1}_{j,n-1}e_{jn}
 +( u_{j1}v_{i1}- v_{j1}u_{i1} )c^{j1}_{i,n-1}e_{in} )$$
$$=\sum_{i=1}^{m} \sum_{j=1,j\ne i}^{m}
( u_{i1}v_{j1}- v_{i1}u_{j1}
)(c^{i1}_{j,n-1}e_{jn}-c^{j1}_{i,n-1}e_{in} )$$
$$=\sum_{i=1}^{m} \sum_{j=1,j\ne i}^{m}
( u_{i1}v_{j1}- v_{i1}u_{j1} ) \left( c^{i1}_{j,n-1}\sum_{k=1}^{r}
b_{jk}e_{kn} -c^{j1}_{i,n-1}\sum_{k=1}^{r} b_{ik}e_{kn} \right)$$
$$=\sum_{i=1}^{m} \sum_{j=1,j\ne i}^{m}
( u_{i1}v_{j1}- v_{i1}u_{j1} ) \left( \sum_{k=1}^{r} (
c^{i1}_{j,n-1}b_{jk} -c^{j1}_{i,n-1} b_{ik} ) e_{kn} \right)=0.
$$

\hfill $\Box$

                          \vspace{6pt}

 By Theorem 2.1, it is easy to know that there exists an
  $h_1\in {\rm Der} N$  such that  the matrix of  $h_1$ relative to
  $\{ e_{10}, e_{11}, \ldots, e_{m0},  e_{m1} \}$ is
  $ {\rm diag} ( -2, n-2,  \ldots, -2m, m(n-2)  ).$
  Applying  Lemma 1.2 we deduce  the following proposition.

 \vspace{6pt}

  \noindent  {\bf  Proposition 2.1 }\quad  $Q_n$-filiform decomposition
  is unique up to isomorphism.

                     \vspace{6pt}

 \noindent {\bf Remark 2.1}\quad A similar argument as in [11] shows that $N(Q_n,m)$ is a
 nilradical of a complete solvable Lie  algebra  (i.e., centerless
 with only inner derivations).

 \vspace{6pt}

  By (2.2) and (2.3), we know that  the matrix of $ \delta $  relative to the
   basis $\Phi$  is a lower triangular matrix,
   so Der$N$ is a solvable Lie algebra.  Let
   $$ T=\{ \delta \in {\rm Der}N \mid
    \delta ( e_{s0}, \ e_{s1} )=( c^{s0}_{s0}e_{s0}, \
  c^{s1}_{s1}e_{s1} ) \},$$

$$\widetilde{N}=\{ \delta \in {\rm Der}N \mid
  \delta e_{s0}= \sum_{i=2}^{n-1} c^{s0}_{si}e_{si} +\sum_{i=1}^{r}c^{s0}_{in}e_{in}, $$ $$
    \delta e_{s1}=\sum_{i=1}^{d-1}c^{s1}_{s,2i}e_{s,2i}
 +\sum_{i=1}^{m}c^{s1}_{i,n-1}e_{i,n-1}+\sum_{i=1}^{r}c^{s1}_{in}e_{in} \}.$$
  Then
${\rm Der}N=T\dot+\widetilde{N}$.
 Obviously $T$ is an abelian subalgebra, $\widetilde{N}$ is a nilpotent ideal.

             \vspace{6pt}

\noindent {\bf Lemma 2.1 }\quad Let $\delta\in {\rm Der}N$, if $N$
is indecomposable (i.e., cannot be decomposed  into a direct  sum
of ideals of $N$), then $\ \lambda_1=\lambda_2=\cdots=\lambda_m$.

             \vspace{6pt}

  The proof of this lemma is similar to the proof of Lemma 2.2 in [10].

  \vspace{6pt}

\noindent   {\bf Lemma 2.2 }\quad  If $N$ is indecomposable, let $
\tau_{i}\in T$ $( 0 \le i \le m )$  such that the matrices of
 $\tau_{0}$, $\tau_{i}$ relative to
  $\{ e_{10}, e_{11}, \ldots, e_{m0}, e_{m1} \}$ are
   $$M_0={\rm diag}\left(  \begin{array}{ccccccc}
    0 & 1 & 0 & 1 & \cdots & 0 & 1
    \end{array}\right),$$
$$M_i={\rm diag}( \  \underbrace{ 0 \  \ 0 \  \ \cdots \  \ 0 \  \ 0  }_{2(i-1)}  \
   -2 \ \ n-2 \ \underbrace{  0 \  \ 0 \  \ \cdots \  \ 0 \  \ 0 }_{2(m-i)} \  ),$$
 respectively,   then  $\{\tau_{0}, \tau_{1}, \ldots, \tau_{m} \}$ is a basis of
  $T$. Hence $\dim T=m+1.$

  \vspace{6pt}

Next we determine the dimension of $\widetilde{N}$.

Observe that $N(Q_n,m,r)$ admits the following $Q_n$-filiform
decomposition:
$$N(Q_n,m,r)=N_1(Q_n,m_1,1)+N_2(Q_n,m_2,1)+\cdots+N_q(Q_n,m_q,1).$$
Let  $s_i=1+\sum_1^{i-1}m_j$,  $\bar{s}_i=s_i+m_i-1$ $(1\le i\le
q)$. Then
 $\{ e_{s_i0}, e_{s_i1}, \ldots, e_{\bar{s}_i0}, e_{\bar{s}_i1} \}$ is a
 minimal system of generators of $N_i(Q_n,m_i,1)$.
 May assume that $e_{s_in}=e_{s_i+1,n}=\cdots=e_{s_i+m_i-1,n}$,
 and   $\{ e_{1n},  e_{s_2n}, \ldots,  e_{s_rn} \}$
  is a   basis of $c^{n-1}N$.

 For any $\delta \in {\rm Der}N$, by (2.5), we have
 $ c^{s1}_{p,n-1}= c^{p1}_{s,n-1}=0$ if $\{ e_{sn}, e_{pn} \}$ is
  linearly independent,  $c^{s1}_{p,n-1}=c^{p1}_{s,n-1}$ if $\{ e_{sn}, e_{pn} \}$
  is linearly dependent (may assume  $e_{sn}=e_{pn}$).
 Obviously $\{ e_{s_in}, e_{s_jn} \}$ $(i\ne j)$ is linearly independent,
 hence  we have
  $$\delta \left(  \begin{array}{c}
  e_{s_i,1} \\ e_{s_i+1,1} \\ \vdots \\  e_{\bar{s}_i,1}
    \end{array}\right)
  =C_{i}^{1} \left(  \begin{array}{c}
  e_{s_i,1} \\ e_{s_i+1,1} \\ \vdots \\  e_{\bar{s}_i,1}
    \end{array}\right)
  +C_{i}^{2}\left(  \begin{array}{c}
  e_{s_i,2} \\ e_{s_i+1,2} \\ \vdots \\  e_{\bar{s}_i,2}
      \end{array}\right)
+C_{i}^{4}\left(  \begin{array}{c}
  e_{s_i,4} \\ e_{s_i+1,4} \\ \vdots \\  e_{\bar{s}_i,4}
      \end{array}\right)
   $$ $$+\cdots    +   C_{i}^{n-3}\left(  \begin{array}{c}
  e_{s_i,n-3} \\ e_{s_i+1,n-3} \\ \vdots \\  e_{\bar{s}_i,n-3}
      \end{array}\right)
    +   C_{i}^{n-1}\left(  \begin{array}{c}
  e_{s_i,n-1} \\ e_{s_i+1,n-1} \\ \vdots \\  e_{\bar{s}_i,n-1}
      \end{array}\right)
    +C_i^n\left(  \begin{array}{c}
  e_{1n} \\ e_{s_2n} \\ \vdots \\ e_{s_rn}
 \end{array}\right) ,$$
where $C^1_i$, $C^2_i, C^4_i,\ldots, C^{n-2}_i$ are diagonal
matrices, $C_i^{n-1}$ is a symmetric $m_i\times m_i$ matrix.

                       \vspace{6pt}

\noindent  {\bf Lemma 2.3 }\quad  Let    $ \delta^l_{ij}$,
 $ \varepsilon_{ik}^l$, $\zeta_{it}^l$, $ \xi_{it}^l\in \widetilde{N}$, such that
   $$   \delta^l_{ij} \left(  \begin{array}{c}
    e_{i1} \\ e_{j1} \\ e_{s0} \\ e_{s1}
    \end{array}\right)
    =\left(  \begin{array}{c}
     e_{j,n-1} \\ e_{i,n-1} \\ 0 \\ 0
    \end{array}\right), \quad
       \varepsilon_{ik}^l \left(  \begin{array}{c}
    e_{i0} \\ e_{i1} \\ e_{s0} \\ e_{s1}
    \end{array}\right)
    =\left(  \begin{array}{c}
   0 \\ e_{s_i,2k}  \\ 0 \\ 0
    \end{array}\right), $$ $$
     \zeta_{it}^l \left(  \begin{array}{c}
    e_{i0} \\ e_{i1} \\ e_{s0} \\ e_{s1}
    \end{array}\right)
    =\left(  \begin{array}{c}
     e_{s_tn} \\ 0 \\ 0 \\ 0
    \end{array}\right),  \quad
        \xi_{it}^l \left(  \begin{array}{c}
    e_{i0} \\ e_{i1} \\ e_{s0} \\ e_{s1}
    \end{array}\right)
    =\left(  \begin{array}{c}
     0 \\  e_{s_tn} \\ 0 \\ 0
    \end{array}\right),$$
 where  $ 1\le l\le q $, $ 2\le k\le d $ $(n=2d+1)$,
  $ s_l \le i<j \le \bar{s}_{l} $,  $ 1\le t\le r. $
      Then for any $\delta\in \widetilde{N}$, we have
   $$ \delta=\sum_{l=1}^q\sum_{i=s_l}^{\bar{s}_{l}}
   \sum_{s_l \le i< j \le \bar{s}_{l}}c^{i1}_{j,n-1}\delta^l_{ij}$$
     $$
   +\sum_{l=1}^q\sum_{i=s_l}^{\bar{s}_l}\left(
  -\sum_{j=1}^{n-2} c^{i0}_{i,j+1}{\rm ad}e_{ij}
    +\sum_{t=1}^{r}c^{i0}_{tn}\zeta_{it}^l
  +c^{i1}_{i2}{\rm ad}e_{i0}
  +\sum_{k=2}^{d}c^{i1}_{i,2k} \varepsilon_{ik}^l
  +\sum_{t=1}^{r}c^{i1}_{tn}\xi_{it}^l \right),   $$
  and
  $$\bigcup_{l=1}^q \{ {\rm ad}e_{iu},  \delta^l_{i}, \delta^l_{ij},
    \varepsilon_{ik}^l, \zeta_{it}^l,  \xi_{it}^l \mid 0\le u\le  n-2,
  2\le k\le  d,   s_l \le i< j \le \bar{s}_{l},  1 \le t \le r  \}$$
   is a basis of $\widetilde{N}$.
   Hence $\dim \widetilde{N}=\sum_{l=1}^{q}((2r+n+d-2)m_l+\frac{1}{2}m_l(m_l-1))$.
  In particular, if $N=N(Q_n,m,1)$, then
   $ \{  {\rm ad}e_{i0},  {\rm ad}e_{i1}, \varepsilon_{ik},
   \delta_{ij},  \zeta_{i1}, \xi_{i1} \mid  2\le k\le  d, 1 \le i< j \le m \}$
     is a  basis of $\widetilde{N}$ , so
 $\dim {\rm Der}N(Q_n,m,1)=(n+d+1)m+1+\frac{1}{2}m(m-1)$.

\section  {Isomorphism theorem  }

In this section applying some properties of root vector
decomposition we obtain isomorphism theorem of  quasi
$Q_n$-filiform Lie algebras.

  \vspace{6pt}

By (1.1) there exists a matrix $A$ such that
   $$\left(  \begin{array}{cc} A & I_{m-r}
\end{array}\right)
\left(  \begin{array}{cccc}
e_{1n}  & e_{2n} & \cdots  & e_{mn}
\end{array}\right)^t
=0.$$

 \noindent {\bf Definition 3.1 }\ \  The matrix $(A \ \ I_{m-r} )$ as above is
 called a related matrix of $N(Q_n,m, r)$ with respect to
 $\{ e_{i0}, e_{i1}\ | \  1\le i \le m  \}$.

                            \vspace{6pt}

 \noindent  {\bf Lemma 3.1 }\quad Let
 $(  A \ \ I_{m-r}  )(  e_{1n} \ e_{2n} \ \cdots \ e_{mn}  )^{t}=0$,
 $B$  a $(m-r)\times m$ matrix and $B(  e_{1n} \  e_{2n} \ \cdots \
e_{mn}  )^{t}=0$. If $(A \ \ I_{m-r})$ is a related matrix and
rank$(B)=m-r$, then there exists
 an invertible matrix  $E$ such that   $EB=( A \ \ I_{m-r} )$.

                \vspace{6pt}

  The proof of this lemma is similar to the proof of Lemma 9 in [9].

\vspace{6pt}

 \noindent {\bf Theorem 3.1 }\quad  If $(  A_i \ \ I_{m-r}  )$ is a related matrix
 of $N_i(Q_n,m,r)$ $(i=1,2)$, then  $N_1$ is isomorphic to $N_2$ if and only if
  there exist   invertible   matrix $E$ and
    monomial matrix  $K$  such that
  $E (A_1  \ \  I_{m-r} )K=( A_2 \ \ I_{m-r} )$.

            \vspace{6pt}

  The proof of this theorem is similar to the proof of  Theorem 3.1 in [11].

\section  {On Aut$\mathbf{N(Q_n,m,r)}$ }

In this section we explicitly determine automorphisms of  quasi
$Q_n$-filiform Lie algebra.

\vspace{6pt}

Let  $\rho$ be a linear transformation of $N(Q_n,m,r)$ such that
$$\rho e_{st}=\sum_{i=1}^{m}\sum_{j=0}^{n-1} b^{st}_{ij}e_{ij}+\sum_{j=1}^{r} b^{st}_{in}e_{in},
   \quad t=0,1,$$
$$ \rho e_{st}=[\rho e_{s0},\rho e_{s,t-1}], \quad 2\leq t \leq n-1, $$ $$
\rho e_{sn}=-[\rho e_{s1}, \rho e_{s,n-1}], \quad 1\leq s \leq
r.$$ Then we have
$$\rho e_{s2}=\sum_{i=1}^{m}\sum_{j=1}^{n-2}(b^{s0}_{i0}b^{s1}_{ij}-b^{s0}_{ij}b^{s1}_{i0})e_{i,j+1}
+\sum_{i=1}^{m}\sum_{j=1}^{n-1}(-1)^jb^{s0}_{ij}b^{s1}_{i,n-j}e_{in},$$
$$\rho e_{st}=\sum_{i=1}^{m}
\sum_{j=1}^{n-t}(b^{s0}_{i0})^{t-2}(b^{s0}_{i0}b^{s1}_{ij}-b^{s0}_{ij}b^{s1}_{i0})e_{i,j+t-1}
$$ $$+\sum_{i=1}^{m}\sum_{j=1}^{n-t+1}(-1)^j
b^{s0}_{ij}(b^{s0}_{i0})^{t-3}(b^{s0}_{i0}b^{s1}_{i,n-j-t+2}-b^{s0}_{i,n-j-t+2}b^{s1}_{i0})e_{in},
\ \ 2<t<n,$$
$$\rho e_{sn}=\sum_{i=1}^{m}b^{s1}_{i1}
(b^{s0}_{i0})^{n-3}(b^{s0}_{i0}b^{s1}_{i1}-b^{s0}_{i1}b^{s1}_{i0})e_{in},
\quad 1\leq s \leq r. $$

                       \vspace{6pt}

\noindent  {\bf Theorem 4.1 }\quad  Let $\rho$ be as above, $(
e_{1n} \ e_{2n} \ \cdots \ e_{mn}  ) =( e_{1n} \ e_{2n} \ \cdots \
e_{rn}  )( I \ \ B )$, then $ \rho \in {\rm Aut} N$ if and only if
the following conditions hold:

(1) For any $s$, there exists only one integer $q_s$ $(1\le q_s\le
m)$  such that

$$\rho e_{s0}=b^{s0}_{q_s0}e_{q_s0}+ \sum_{i=2}^{n-1}b^{s0}_{q_si}e_{q_si}
+\sum_{i=1}^{r}b^{s0}_{in}e_{in},
$$
$$\rho e_{s1}= \sum_{i=1}^{n-2} b^{s1}_{q_si}e_{q_si}
+\sum_{i=1}^{m}b^{s1}_{i,n-1}e_{i,n-1}+\sum_{i=1}^{r}b^{s1}_{in}e_{in}.
$$

(2) There exists a permutation  matrix $T$ such that
$$
\left(  \begin{array}{cccc}  q_1 & q_2 & \cdots & q_m
\end{array}\right)=
\left(  \begin{array}{cccc} 1 & 2 &  \cdots & m
\end{array}\right)T.$$

(3) For any $ \ 1\le s, \ p\le m$,
$$ b^{s0}_{q_s0}b^{s1}_{q_s1}\ne 0, $$ $$
b^{s1}_{q_s1}b^{p1}_{q_s,n-1}e_{q_sn}=b^{s1}_{q_p,n-1}b^{p1}_{q_p1}e_{q_pn}.$$

(4)
$$\sum_{j=1}^{p}(-1)^jb^{s1}_{q_sj}b^{s1}_{q_s,p-j+1}=0, \quad p=3, 5, \ldots, n-2.$$

(5)
$$\left(  \begin{array}{cc} I & B
\end{array}\right)T_2K_2=
\left(  \begin{array}{cc} I & B
\end{array}\right) T_1K_1B.$$
where $T=( T_1 \ \ T_2 )$, $T_1$ is a $m\times r$ matrix,
$K_1={\rm diag}( k_1, \ldots, k_r )$, $K_2={\rm diag}( k_{r+1},
\ldots, k_m )$, $k_i=(b^{i0}_{q_i0})^{n-2}(b^{i1}_{q_i1})^2$.

         \vspace{6pt}

\noindent  {\bf Proof. }\quad $(\Rightarrow)$: As  $\rho \in {\rm
Aut} N$, for any $s$,
$$\rho e_{sn}=\sum_{i=1}^{m}b^{s1}_{i1}
(b^{s0}_{i0})^{n-3}(b^{s0}_{i0}b^{s1}_{i1}-b^{s0}_{i1}b^{s1}_{i0})e_{in}.$$
Note that $ \rho e_{sn}\ne 0$, there exists an integer $q_s$ such
that
$$b^{s0}_{q_s0}b^{s1}_{q_s1}-b^{s1}_{q_s0}b^{s0}_{q_s1}\ne 0.  \eqno (4.1)$$

For any $ \ 1\le s \ne p\le m$, $ 0\le t, k \le 1$, by $\rho
[e_{st}, e_{pk}]=0$, we have
  $$\sum_{i=1}^{m} \sum_{j=1}^{n-2}
 (b^{st}_{i0}b^{pk}_{ij}-b^{st}_{ij}b^{pk}_{i0})e_{i,j+1}+
\sum_{i=1}^{m}\sum_{j=1}^{n-1} (-1)^jb^{st}_{ij}
b^{pk}_{i,n-j}e_{in}=0.
 \eqno (4.2)$$

Now we prove that for any $1\le p\ne s\le m$,
$$b^{pt}_{q_st}=b^{p,1-t}_{q_st}=0, \quad  t=0,1.  \eqno (4.3)$$

 By (4.2), we have
$b^{s0}_{i0}b^{p0}_{i1}-b^{s0}_{i1}b^{p0}_{i0}
=b^{s1}_{i0}b^{p1}_{i1}-b^{s1}_{i1}b^{p1}_{i0}=0$, $\forall i$,
then for $t=0, 1$,
$$\rho( b^{st}_{q_st}e_{pt}-b^{pt}_{q_st}e_{st} )
-\sum_{i=1}^{r}(
b^{st}_{q_st}b^{pt}_{in}-b^{pt}_{q_st}b^{st}_{in})e_{in}
$$ $$=b^{st}_{q_st}\left( \sum_{j=0}^{n-1}
 b^{pt}_{q_sj}e_{q_sj}
+\sum_{i=1,i\ne q_s}^{m} \sum_{j=0}^{n-1}
 b^{pt}_{ij}e_{ij} \right)
-b^{pt}_{q_st}\left( \sum_{j=0}^{n-1}
 b^{st}_{q_sj}e_{q_sj}
+\sum_{i=1,i\ne q_s}^{m} \sum_{j=0}^{n-1}
 b^{st}_{ij}e_{ij} \right)
$$  $$
= \sum_{j=2}^{n-1}(
 b^{st}_{q_st}b^{pt}_{q_sj}
-b^{pt}_{q_st} b^{st}_{q_sj} )e_{q_sj} + \sum_{i=1,i\ne q_s}^{m}
\sum_{j=0}^{n-1} ( b^{st}_{q_st}b^{pt}_{ij}
-b^{pt}_{q_st}b^{st}_{ij} )e_{ij},
$$
therefore
$$
[\rho(b^{st}_{q_st}e_{pt}-b^{pt}_{q_st}e_{st}), \rho e_{s,1-t}]$$
$$ = -\sum_{j=2}^{n-2}(
 b^{st}_{q_st}b^{pt}_{q_sj}-b^{pt}_{q_st} b^{st}_{q_sj} )b^{s,1-t}_{q_s0}e_{q_s,j+1}
 +\sum_{j=2}^{n-1}(-1)^j(
 b^{st}_{q_st}b^{pt}_{q_sj}-b^{pt}_{q_st} b^{st}_{q_sj} )b^{s,1-t}_{q_s,n-j}e_{q_sn}
$$ $$-
\sum_{i=1,i\ne q_s}^{m} \sum_{j=1}^{n-2}
(b^{st}_{q_st}b^{pt}_{ij}-b^{pt}_{q_st}b^{st}_{ij})b^{s,1-t}_{i0}e_{i,j+1}
-\sum_{i=1,i\ne q_s}^{m}  \sum_{j=1}^{n-2} (
b^{st}_{q_st}b^{pt}_{i0} -b^{pt}_{q_st}b^{st}_{i0}
)b^{s,1-t}_{ij}e_{i,j+1}
$$
$$+\sum_{i=1,i\ne q_s}^{m}\sum_{j=1}^{n-1}(-1)^j ( b^{st}_{q_st}b^{pt}_{i0}
-b^{pt}_{q_st}b^{st}_{i0} )b^{s,1-t}_{i,n-j}e_{in}. $$
But on the
other hand, for $t=0,1$,
$$[\rho(b^{st}_{q_st}e_{pt}-b^{pt}_{q_st}e_{st}), \rho e_{s,1-t}]
=\pm b^{pt}_{q_st}\rho e_{s2} $$
 $$ =\pm b^{pt}_{q_st}
\sum_{i=1}^{m}
\left(\sum_{j=1}^{n-2}(b^{s0}_{i0}b^{s1}_{ij}-b^{s0}_{ij}b^{s1}_{i0})e_{i,j+1}
+\sum_{j=1}^{n-1}(-1)^jb^{s0}_{ij}b^{s1}_{i,n-j}e_{in} \right). $$
 Comparing the coefficients of $e_{q_s2}$, we have
$$b^{pt}_{q_st}
(b^{s0}_{q_s0}b^{s1}_{q_s1}-b^{s1}_{q_s0}b^{s0}_{q_s1})=0.$$
 By (4.1), we have
$$b^{pt}_{q_st}=0, \quad p\ne s, \ t=0,1.$$

Similarly by
$$[\rho( b^{st}_{q_st}e_{p,1-t}-b^{p,1-t}_{q_st}e_{st}), \rho e_{s,1-t}]
=\mp b^{p,1-t}_{q_st}\rho e_{s2}, \quad t=0,1,$$
 we have
$$b^{p,1-t}_{q_st}=0, \quad p\ne s, \ t=0,1.$$
Thus (4.3) holds.

If  $\exists s\ne p$ such that $q_s=q_p$, by (4.3),
 $b^{i0}_{q_s0}=b^{i1}_{q_s0}=0$, $\forall i$.   This implies that $e_{q_s0}\not\in \rho (N)$,
 a contradiction. So there exits a permutation  matrix $T$ such that
 $$
\left(  \begin{array}{cccc}  q_1 & q_2 & \cdots & q_m
\end{array}\right)=
\left(  \begin{array}{cccc} 1 & 2 &  \cdots & m
\end{array}\right)T. \eqno (4.4)$$

By (4.3), (4.4), and $\rho [e_{s1}, e_{s2}]=0$, we have
$$ 0= \sum_{j=1}^{n-3} b^{s1}_{q_s0}
(b^{s0}_{q_s0}b^{s1}_{q_sj}-b^{s1}_{q_s0}b^{s0}_{q_sj})e_{q_s,j+2}
$$ $$+\sum_{j=1}^{n-2} (-1)^j b^{s1}_{q_sj}
(b^{s0}_{q_s0}b^{s1}_{q_s,n-j-1}-b^{s1}_{q_s0}b^{s0}_{q_s,n-j-1})e_{q_sn}.$$
Hence $b^{s1}_{q_s0}=0$, then
$$\rho e_{sn}
=(b^{s0}_{q_s0})^{n-2}(b^{s1}_{q_s1})^2e_{q_sn},$$
 so
$$b^{s0}_{q_s0}b^{s1}_{q_s1}\ne 0.$$

By
$$0=\rho
[e_0,e_{n-1}]=-b^{s0}_{q_s1}(b^{s0}_{q_s0})^{n-2}b^{s1}_{q_s1}e_{q_s,n},$$
we have  $$b^{s0}_{q_s1}=0.$$

For $t=2, 4, \ldots, n-3$, by $$0=\rho [e_{s1},
e_{st}]=(b^{s0}_{q_s0})^{t-1}\sum_{j=1}^{n-t}(-1)^j
b^{s1}_{q_sj}b^{s1}_{q_s,n-t-j+1}e_{q_s,n},$$ we know that (4)
holds.

 By (4.2) and (4.3), when $t=0$, $k=1$, we have
$$ \sum_{j=1}^{n-2}
 b^{s0}_{q_s0}b^{p1}_{q_sj}e_{q_s,j+1}+
\sum_{i=1}^{m}\sum_{j=1}^{n-1} (-1)^jb^{s0}_{ij}
b^{p1}_{i,n-j}e_{in}=0, $$ then
$$b^{p1}_{q_sj}=0, \quad 1\le j\le n-2.$$
 Hence
 $$b^{s0}_{q_p,n-1} b^{p1}_{q_p1}e_{q_pn}=0,$$
  so
$$b^{s0}_{q_p,n-1}=0.$$

By (4.2) and (4.3), when $t=k=0$, we have
$$b^{s0}_{q_pj}=b^{p0}_{q_sj}=0, \quad 1\le j\le n-2.$$

By (4.2) and (4.3), when $t=k=1$, we have
$$b^{s1}_{q_s1}b^{p1}_{q_s,n-1}e_{q_sn}= b^{s1}_{q_p,n-1}b^{p1}_{q_p1}e_{q_pn}.$$

Now we have showed that (1), (2), (3) and (4) hold.

At last  we show that (5) holds. Let
$$ K={\rm diag}(k_1, \ k_2, \ \ldots, \ k_m ),
\quad k_s=(b^{s0}_{q_s0})^{n-2}(b^{s1}_{q_s1})^{2}.$$
Then
$$
\rho \left(  \begin{array}{cccc}  e_{1n} & e_{2n} & \cdots &
e_{mn}
\end{array}\right)=
\left(  \begin{array}{cccc}  e_{1n} & e_{2n} & \cdots & e_{mn}
\end{array}\right)TK$$
$$=\left(  \begin{array}{cccc}  e_{1n} & e_{2n} & \cdots & e_{rn}
\end{array}\right)
\left(  \begin{array}{cc} I & B
\end{array}\right)TK.
$$

 But on the other hand,
$$\rho \left(  \begin{array}{cccc}  e_{1n} & e_{2n} & \cdots &
e_{mn}
\end{array}\right)
=\rho \left(  \begin{array}{cccc}  e_{1n} & e_{2n} & \cdots &
e_{rn}
\end{array}\right)
\left(  \begin{array}{cc} I & B
\end{array}\right)$$
$$=\left(  \begin{array}{cccc} e_{q_1n} & e_{q_2n} & \cdots & e_{q_rn}
\end{array}\right)K_1
\left(  \begin{array}{cc} I & B
\end{array}\right)
$$ $$=
\left(  \begin{array}{cccc}  e_{1n} & e_{2n} & \cdots & e_{mn}
\end{array}\right)T_1K_1
\left(  \begin{array}{cccc} I & B
\end{array}\right)
$$
$$=\left(  \begin{array}{cccc}  e_{1n} & e_{2n} & \cdots & e_{rn}
\end{array}\right)
\left(  \begin{array}{cc} I & B
\end{array}\right)
T_1K_1 \left(  \begin{array}{cc} I & B
\end{array}\right)$$
where $T=( T_1 \ \ T_2 )$, $T_1$ is an $m\times r$ matrix,
$K_1={\rm diag}( k_1, \ldots, k_r )$. So we have
$$\left(  \begin{array}{cc} I & B
\end{array}\right)TK
=\left(  \begin{array}{cc} I & B
\end{array}\right)
T_1K_1 \left(  \begin{array}{cc} I & B
\end{array}\right).
$$
Let  $K_2={\rm diag}( k_{r+1}, \ldots, k_m )$, then we have
$$\left(  \begin{array}{cc} I & B
\end{array}\right)T_2K_2
=\left(  \begin{array}{cc} I & B
\end{array}\right)
T_1K_1B.
$$

    \vspace{6pt}

 $(\Leftarrow)$: We only prove that $\rho [x, y]=[\rho x,  \rho y]$ for any $x, y \in N$.
 Set
   $$x=\sum_{i=1}^{m} \sum_{j=0}^{n-1} u_{ij}e_{ij}
   +\sum_{i=1}^{r} u_{in}e_{in}, \quad \quad
  y=\sum_{i=1}^{m}\sum_{j=0}^{n-1} v_{ij}e_{ij}
     +\sum_{i=1}^{r} v_{in}e_{in}.  $$

For any $s$, by (5), we have
$$\rho e_{sn} =(b^{s1}_{q_s1})^2
(b^{s0}_{q_s0})^{n-2}e_{q_sn}.$$

If $\rho [e_{st},  e_{pk}]=[\rho e_{st}, \rho e_{pk}]$ for any
$e_{st}, e_{pk}$,  then we have
$$\rho [x, y]=\rho \sum_{i=1}^{m} \left[  \sum_{j=0}^{n-1} u_{ij}e_{ij}, \
       \sum_{j=0}^{n-1} v_{ij}e_{ij} \right] $$ $$
=\rho \sum_{i=1}^{m}\left(u_{i0}\sum_{j=1}^{n-2}v_{ij}e_{i,j+1}
+\sum_{j=1}^{n-1}(-1)^ju_{ij}v_{i,n-j}e_{in}
-v_{i0}\sum_{j=1}^{n-2} u_{ij}e_{i,j+1} \right)
$$ $$=\sum_{i=1}^{m}\left(\sum_{j=1}^{n-2}(u_{i0}v_{ij}-v_{i0}u_{ij})\rho e_{i,j+1}
+\sum_{j=1}^{n-1}(-1)^ju_{ij}v_{i,n-j}\rho e_{in}\right)$$
$$=\left[\sum_{i=1}^{m} \sum_{j=0}^{n-1} u_{ij}\rho
e_{ij}, \ \sum_{i=1}^{m}\sum_{j=0}^{n-1} v_{ij}\rho
e_{ij}\right]$$
$$ =[\rho x, \rho y].$$

 We now  prove that
$\rho [e_{st},  e_{pk}]=[\rho e_{st}, \rho e_{pk}]$ for any
$e_{st}, e_{pk}$.

 Obviously $\rho [e_{st},  e_{pk}]=0=[\rho e_{st},
\rho e_{pk}]$ if  $t=n$ or $k=n$.

Obviously $\rho [e_{st}, e_{pk}]=0=[\rho e_{st}, \rho e_{pk}]$ if
$s\ne p$, and $t\neq 1$ or $k\neq1$.

If $s\ne p$, $t=k=1$, by (3), we have
  $\rho [e_{st},  e_{pk}]=0=[\rho e_{st},\rho e_{pk}]$.

Next we  prove that
$$\rho [e_{st},e_{sk}]=[\rho e_{st}, \rho e_{sk}], \quad  0\leq t < k <n.$$

Case 1:  $t=0$. This equation is obvious.

Case 2: $t=1$. By (4) and note that
$\sum_{j=1}^{p}(-1)^jb^{s1}_{q_sj}b^{s1}_{q_s,p-j+1}=0$ when $p$
is even, this equation holds.

Case 3: $t>d$ $(n=2d+1)$. This equation is obvious.

Case 4: $t=d$.
$$[\rho e_{sd}, \rho e_{sk}]=0=\rho [e_{sd},e_{sk}],
\quad k > d+1,$$
$$[\rho e_{sd}, \rho
e_{s,d+1}]=(-1)^{d}(b^{s0}_{q_s0})^{n-2}(b^{s1}_{q_s1})^2e_{q_sn}=\rho
[e_{sd},e_{s,d+1}].$$

Case 5: $1<t<d$ and $k<d$.  Obviously $\rho[e_{st},e_{sk}]=0$.

Set $t=d-t'$, $k=d-k'$. By (4) and note that
$\sum_{j=1}^{p}(-1)^jb^{s1}_{q_sj}b^{s1}_{q_s,p-j+1}=0$ when $p$
is even,  we have
$$[\rho e_{s,d-t'}, \rho e_{s,d-k'}]$$ $$
=\left[\sum_{j=1}^{n-d+t'}(b^{s0}_{q_s0})^{d-t'-1}b^{s1}_{q_sj}e_{q_s,j+d-t'-1},
\sum_{j=1}^{n-d+k'}(b^{s0}_{q_s0})^{d-k'-1}b^{s1}_{q_sj}e_{q_s,j+d-k'-1}\right]$$
$$=(b^{s0}_{q_s0})^{n-k'-t'-3} \left[\sum_{j=1}^{t'+1}b^{s1}_{q_sj}e_{q_s,j+d-t'-1},
\sum_{j=k'+2}^{k'+t'+2}b^{s1}_{q_sj}e_{q_s,j+d-k'-1} \right]$$
$$+(b^{s0}_{q_s0})^{n-k'-t'-3} \left[\sum_{j=t'+2}^{k'+t'+2}b^{s1}_{q_sj}e_{q_s,j+d-t'-1},
\sum_{j=1}^{k'+1}b^{s1}_{q_sj}e_{q_s,j+d-k'-1} \right]$$
$$=(b^{s0}_{q_s0})^{n-k'-t'-3}\sum_{j=1}^{k'+t'+2}(-1)^jb^{s1}_{q_sj}b^{s1}_{q_s,k'+t'+2-j+1}e_{q_sn}
=0.$$

Case 6: $1<t<d$ and $k\geq d$. Obviously $\rho[e_{st},e_{sk}]=0$.
Set $t=d-t'$, $k=d+k'$.

If $t'+1\geq k'$,
$$[\rho e_{s,d-t'}, \rho e_{s,d+k'}]$$ $$
=\left[\sum_{j=1}^{n-d+t'}(b^{s0}_{q_s0})^{d-t'-1}b^{s1}_{q_sj}e_{q_s,j+d-t'-1},
\sum_{j=1}^{n-d-k'}(b^{s0}_{q_s0})^{d+k'-1}b^{s1}_{q_sj}e_{q_s,j+d+k'-1}\right]$$
$$=\left[\sum_{j=1}^{t'-k'+2}(b^{s0}_{q_s0})^{d-t'-1}b^{s1}_{q_sj}e_{q_s,j+d-t'-1},
\sum_{j=1}^{t'-k'+2}(b^{s0}_{q_s0})^{d+k'-1}b^{s1}_{q_sj}e_{q_s,j+d+k'-1}\right]$$
$$=(b^{s0}_{q_s0})^{n-t'+k'-3}\sum_{j=1}^{t'-k'+2}(-1)^jb^{s1}_{q_sj}b^{s1}_{q_s,t'-k'+2-j+1}e_{q_sn}
=0.$$

If $t'+1< k'$,
$$[\rho e_{s,d-t'}, \rho e_{s,d+k'}]$$ $$
=\left[\sum_{j=1}^{n-d+t'}(b^{s0}_{q_s0})^{d-t'-1}b^{s1}_{q_sj}e_{q_s,j+d-t'-1},
\sum_{j=1}^{n-d-k'}(b^{s0}_{q_s0})^{d+k'-1}b^{s1}_{q_sj}e_{q_s,j+d+k'-1}\right]=0.
$$

 \hfill $\Box$

   \end{document}